\documentclass{amsart}

\usepackage[leqno]{amsmath}
\usepackage{latexsym,amsthm,amsfonts,amssymb,hyperref}
\usepackage[all]{xy} \SelectTips{eu}{}

\newcommand{\numberseries}{\mdseries}   
\newlength{\thmtopspace}                
\newlength{\thmbotspace}                
\newlength{\thmheadspace}               
\newlength{\thmindent}                  
\newtheoremstyle{bfupright head,upright body}
                {\thmtopspace}{\thmbotspace}
                {\upshape}{\thmindent}{\bfseries}{.}{\thmheadspace}
                {{\numberseries \thmnumber{(#2) }}\thmnote{#3}}

\newtheoremstyle{fixed bf head,slanted body}
                {\thmtopspace}{\thmbotspace}{\slshape}
                {\thmindent}{\bfseries}{.}{\thmheadspace}
                {{\numberseries \thmnumber{(#2) }}\thmname{#1}\thmnote{ (#3)}}

\newtheoremstyle{fixed bf head,upright body}
                {\thmtopspace}{\thmbotspace}{\upshape}
                {\thmindent}{\bfseries}{.}{\thmheadspace}
                {{\numberseries \thmnumber{(#2) }}\thmname{#1}\thmnote{ (#3)}}

\newtheoremstyle{numbered paragraph}
                {\thmtopspace}{\thmbotspace}{\upshape}
                {\thmindent}{\upshape}{}{0pt}
                {{\numberseries \thmnumber{(#2) }}}

\theoremstyle{bfupright head,upright body}
\newtheorem{res}{}[section]             \newtheorem*{res*}{}
\newtheorem{bfhpg}[res]{}               \newtheorem*{bfhpg*}{}

\theoremstyle{fixed bf head,slanted body}
\newtheorem{thm}[res]{Theorem}          \newtheorem*{thm*}{Theorem}
\newtheorem{prp}[res]{Proposition}      \newtheorem*{prp*}{Proposition}
\newtheorem{lem}[res]{Lemma}            \newtheorem*{lem*}{Lemma}

\theoremstyle{fixed bf head,upright body}
       \newtheorem*{dfn*}{Definition}
           \newtheorem*{rmk*}{Remark}

\theoremstyle{numbered paragraph}
\newtheorem{ipg}[res]{}

\newlength{\thmlistleft}        
\newlength{\thmlistright}       
\newlength{\thmlistpartopsep}   
\newlength{\thmlisttopsep}      
\newlength{\thmlistparsep}      
\newlength{\thmlistitemsep}     

\newcounter{prt}
\newenvironment{prt}{\begin{list}{\upshape (\alph{prt})}%
    {\usecounter{prt}%
      \setlength{\leftmargin}{\thmlistleft}%
      \setlength{\labelwidth}{\thmlistleft}%
      \setlength{\rightmargin}{\thmlistright}%
      \setlength{\partopsep}{\thmlistpartopsep}%
      \setlength{\topsep}{\thmlisttopsep}%
      \setlength{\parsep}{\thmlistparsep}%
      \setlength{\itemsep}{\thmlistitemsep}}}%
  {\end{list}}%

\newenvironment{prf*}[1][Proof]{%
  \begin{proof}[\bf #1]
    \setcounter{equation}{0}
    }
  {\end{proof}
}
\renewcommand{\eqref}[1]{\pgref{eq:#1}}
\newcommand{\pgref}[1]{(\ref{#1})}
\newcommand{\thmref}[2][Theorem~]{#1\pgref{thm:#2}}
\newcommand{\prpref}[2][Proposition~]{#1\pgref{prp:#2}}
\newcommand{\lemref}[2][Lemma~]{#1\pgref{lem:#2}}
\newcommand{\thmcite}[2][?]{\cite[thm.~#1]{#2}}
\newcommand{\prpcite}[2][?]{\cite[prop.~#1]{#2}}

\def\urltilda{\kern -.15em\lower .7ex\hbox{\~{}}\kern .04em} 
\newcommand{\set}[2][\mspace{1mu}]{\{#1 #2 #1\}}
\newcommand{\setof}[3][\mspace{1mu}]{\{#1#2 \mid #3#1\}}
\newcommand{\kk}{\Bbbk}
\newcommand{\qtext}[1]{\quad\text{#1}\quad}
\newcommand{\qqtext}[1]{\qquad\text{#1}\qquad}
\newcommand{\qand}{\qtext{and}}
\newcommand{\qqand}{\qqtext{and}}
\newcommand{\deq}{\:=\:}
\newcommand{\dle}{\:\le\:}
\renewcommand{\a}{\alpha}
\renewcommand{\b}{\beta}
\newcommand{\n}{\mathfrak{n}}
\newcommand{\is}{\cong}
\newcommand{\lra}{\longrightarrow}
\newcommand{\xra}[2][]{\xrightarrow[#1]{\;#2\;}}
\newcommand{\pows}[2][k]{#1[\mspace{-2.3mu}[#2]\mspace{-2.3mu}]}
\newcommand{\dmapdef}[4][\lra]{\nobreak{#2\colon #3\:#1\:#4}}
\newcommand{\cls}[1]{[#1]}
\renewcommand{\H}[2][]{\operatorname{H}_{#1}(#2)}

\newcommand{\rnk}[2][k]{\operatorname{rank}_{#1}#2}
\newcommand{\type}[2][R]{\operatorname{type}_{#1}#2}
\newcommand{\Hom}[3][R]{\operatorname{Hom}_{#1}(#2,#3)}
\newcommand{\tp}[3][R]{\nobreak{#2\otimes_{#1}#3}}
\newcommand{\Tor}[4][R]{\operatorname{Tor}^{#1}_{#2}(#3,#4)}
\newcommand{\clT}{\mathbf{T}}
\newcommand{\clB}{\mathbf{B}}
\newcommand{\clG}[1]{\mathbf{G}(#1)}
\newcommand{\clH}[1]{\mathbf{H}(#1)}
\newcommand{\Hilb}[1]{\operatorname{Hilb}_{#1}(t)}
\newcommand{\pf}[2]{\operatorname{pf}_{#1}(#2)}
\newcommand{\e}[1]{\varepsilon_{#1}}
\newcommand{\Ox}{O_x} 
\newcommand{\Oy}{{}^yO}
\newcommand{\fa}{\mathfrak{a}}
\newcommand{\fb}{\mathfrak{b}}
\newcommand{\fg}{\mathfrak{g}}
\numberwithin{equation}{res}

\begin{document}

\title{Trimming a Gorenstein ideal}

\author[L.\,W. Christensen]{Lars Winther Christensen}

\address{Texas Tech University, Lubbock, TX 79409, U.S.A.}

\email{lars.w.christensen@ttu.edu}

\urladdr{http://www.math.ttu.edu/\urltilda lchriste}

\author[O. Veliche]{Oana Veliche}

\address{Northeastern University, Boston, MA~02115, U.S.A.}

\email{o.veliche@neu.edu}

\author[J. Weyman]{Jerzy Weyman}

\address{University of Connecticut, Storrs, CT~06269, U.S.A.}

\email{jerzy.weyman@uconn.edu}

\urladdr{http://www.math.uconn.edu/\urltilda weyman}

\thanks{This work is part of a body of research that started during
  the authors' visit to MSRI in Spring 2013 and continued during a
  months-long visit by L.W.C.\ to Northeastern University; the
  hospitality of both institutions is acknowledged with
  gratitude. L.W.C.\ was partly supported by NSA grant H98230-11-0214,
  and J.W.\ was partly supported by NSF DMS grant 1400740.}

\date{19 January 2017}



\begin{abstract}
  Let $Q$ be a regular local ring of dimension $3$. We show how to
  trim a Gorenstein ideal in $Q$ to obtain an ideal that defines a
  quotient ring that is close to Gorenstein in the sense that its
  Koszul homology algebra is a Poincar\'e duality algebra $P$ padded
  with a non-zero graded vector space on which $P_{\ge 1}$ acts
  trivially. We explicitly construct an infinite family of such rings.
\end{abstract}

\maketitle

\thispagestyle{empty}

\section{Introduction}
\label{sec:introduction}

\noindent
Let $Q$ be a regular local ring with maximal ideal $\n$. Quotient
rings of $Q$ that have projective dimension at most $3$ as $Q$-modules
have been classified based on the multiplicative structure of their
Koszul homology algebras. To be precise, let $\fa \subseteq \n^2$ be
an ideal such that the minimal free resolution of $R=Q/\fa$ over $Q$
has length at most $3$. By a result of Buchsbaum and Eisenbud
\cite{DABDEs77}, the resolution carries a structure of an associative
differential graded commutative algebra, and based on that structure
Avramov, Kustin, and Miller \cite{AKM-88} and Weyman \cite{JWm89}
established a classification in terms of the induced multiplicative
structure on $\Tor[Q]{*}{R}{\kk}$, where $\kk$ is the residue field of
$Q$. Finally, as graded $\kk$-algebras, the Koszul homology algebra of
$R$ and $\Tor[Q]{*}{R}{\kk}$ are isomorphic; see Avramov~\cite{LLA12}
for an in-depth treatment.

An ideal $\fa \subset Q$ is called \emph{Gorenstein} if the quotient
$R=Q/\fa$ is a Gorenstein ring. By a classic result of Avramov and
Golod \cite{LLAESG71}, a Gorenstein ring is characterized by the fact
that its Koszul homology algebra $A=\H{K^R}$ has Poincar\'e
duality. In the classification mentioned above, a Gorenstein ring that
is not complete intersection belongs to a parametrized family
$\clG{r}$, where $r$ is the rank of the canonical map
\begin{equation*}
  \dmapdef{\delta}{A_2}{\Hom[\kk]{A_1}{A_3}}\:;
\end{equation*}
see \cite[1.4.2]{LLA12}. It was conjectured in \cite{LLA12} that all
rings of class $\clG{r}$ are Gorenstein, but Christensen and Veliche
\cite{LWCOVl14} gave sporadic examples of rings of class $\clG{r}$
that are not Gorenstein. In this paper we present a systematic
construction and achieve:

\begin{thm}
  \label{thm:a}
  Let $Q$ be the power series algebra in three variables over a
  field. For every $r \ge 3$ there is quotient ring of $Q$ that is of
  class $\clG{r}$ and not~Gorenstein.
\end{thm}

The quotient rings in \thmref{a} are obtained as follows: Let $\n$ be
the maximal ideal of $Q$ and start with a graded Gorenstein ideal $\fg
\subseteq \n^2$ generated by $2m+1$ elements.  Trim $\fg$ by replacing
one minimal generator $g$ by $\n g$; this removes a 1-dimensional
subspace from $\fg$. The quotient of $Q$ by the resulting ideal is a
ring of type $2$; in particular, it is not Gorenstein, and for $m \ge
3$ it is of class $\clG{r}$. \thmref{a} is consequence of \prpref{a},
which builds on a more general but slightly less precise statement
about local rings, \thmref{2}.

\section{Local rings}
\label{sec:G from Gor}

\noindent
Let $Q$ be a $d$-dimensional regular local ring with maximal ideal
$\n$ and residue field $\kk$.  For an ideal $\fa$ in $Q$, we denote by
$\mu(\fa)$ the minimal number of generators of $\fa$. Let $\fa
\subseteq \n^2$ be an ideal and set $R=Q/\fa$. We denote by $K^R$ the
Koszul complex on a minimal set of generators for the maximal ideal
$\n/\fa$ of $R$; one has $K^R = \tp[Q]{R}{K^Q}$. The Koszul complex is
an exterior algebra, and the homology algebra $A=\H{K^R}$ is a
graded-commutative $\kk$-algebra. Denote by $c$ the projective
dimension of $R$ as a $Q$-module; by the Auslander--Buchsbaum Formula
and depth sensitivity of the Koszul complex one has $c =
\max\setof{i}{A_i\ne 0}$.  The number $\rnk[\kk]{(A_c)}$ is called the
\emph{type} of $R$. If the ideal $\fa$ is $\n$-primary, then one has
$c=d$ and the type of $R$ is the socle rank, i.e.\ $\type[]{(R)} =
\rnk[\kk]{(0 :_R \n/\fa)}$.

\begin{bfhpg}[Classification]
  \label{class}
  Let $Q$ be as above, and let $\fa \subseteq \n^2$ be an ideal such
  that $R = Q/\fa$ has projective dimension $3$ as a $Q$-module. The
  possible multiplicative structures on the graded-commutative
  $\kk$-algebra $A = \H{K^R} \is \Tor[Q]{*}{R}{\kk}$ were identified
  in \cite{AKM-88}.  By assumption one has $A_{\ge 4} =0$, and the
  possible structures are described by the invariants
  \begin{equation*}
    p = \rnk[\kk]{(A_1\cdot A_1)}\,, \ \ 
    q = \rnk[\kk]{(A_1\cdot A_2)}\,, \ \ \text{and}\, \ \   
    r = \rnk[\kk]{(A_2 \xra{\delta} \Hom[\kk]{A_1}{A_3})}\,.
  \end{equation*}
  From \thmcite[3.1]{LLA12} one extracts the following description of
  all the possible classes of rings that are not Gorenstein.
  \begin{equation}
    \label{eq:tab}
    \begin{array}{r|ccc|l}
      \text{Class} &  p & q & r & \text{Restrictions}\\
      \hline
      \clB & 1 &1 &2 \\
      \clG{r} & 0 &1 &r & 2 \le r \le \mu(\fa)-2\\
      \clH{p,q} &  p & q & q & q \le \type[]{(R)}\\
      \clT & 3 &0 &0 \\
    \end{array}
  \end{equation}  
\end{bfhpg}
\noindent
In \cite{AKM-88} the multiplication tables for the different
structures are given. In particular, if $R = Q/\fa$ is a ring of class
$\clG{r}$, then with $m = \mu(\fa)$ and $t = \type[]{(R)}$ there exist
bases for $A_1$, $A_2$, and $A_3$:
\begin{equation*}
  \mathbf{e}_1,\ldots,\mathbf{e}_m\,,\quad \mathbf{f}_1,\ldots,\mathbf{f}_{m+t-1}\,,
  \quad \text{and}\quad \mathbf{g}_1,
  \ldots,\mathbf{g}_t
\end{equation*}
such that the only non-zero products are $\mathbf{e}_i\mathbf{f}_i =
\mathbf{g}_1 = -\mathbf{f}_i\mathbf{e}_i$ for $1\le i \le r$. That is,
the subalgebra $P$ of $A$ spanned by $1$,
$\mathbf{e}_1,\ldots,\mathbf{e}_r$,
$\mathbf{f}_1,\ldots,\mathbf{f}_r$, and $\mathbf{g}_1$ is a pure
Poincar\'e duality algebra, in the sense that the only non-trivial
products are those from the perfect pairing. Moreover, $P_{\ge 1}$
acts trivially on the rest of $A$.

The next result is proved in \cite{CVW-2}; the argument is based on
linkage theory and cannot be reproduced here without significant
overhead.

\begin{prp}
  \label{prp:mu5}
  Let $(Q,\n)$ be a regular local ring and let $\fa \subseteq \n^2$ be
  a perfect ideal of grade $3$ that is minimally generated by $5$
  elements and not Gorenstein.  If, with the notation above, the ring
  $Q/\fa$ has $p=0$, then it has $r \le 1$. \qed
\end{prp}

\begin{lem}
  \label{lem:type}
  Let $(Q,\n)$ be a regular local ring and consider an $\n$-primary
  ideal $\fg \subseteq \n^2$, minimally generated by elements
  $g_0,\ldots,g_k$. Let $s_1,\ldots,s_t$ be elements of $Q$ whose
  classes in $Q/\fg$ form a basis for the socle. The ideal $\fa = \n
  g_0 + (g_1,\ldots,g_k)$ is $\n$-primary, and if $\n s_i \subseteq
  \fa$ holds for all $i=1,\ldots,t$, then the classes of
  $g_0,s_1,\ldots,s_t$ in $Q/\fa$ form a basis for the socle; in
  particular one has $\type[]{(Q/\fa)} = \type[]{(Q/\fg)} +1$.
\end{lem}

\begin{prf*}
  As $\fg$ is $\n$-primary, it follows from the containment $\n
  \fg\subseteq \fa$ that $\fa$ is $\n$-primary. Consider the rings
  $R=Q/\fa$ and $S=Q/\fg$; there is an exact sequence
  \begin{equation*}
    0 \lra \fg/\fa \lra R \lra S \lra 0\:,
  \end{equation*}
  and an isomorphism of $Q$-modules $\fg/\fa\is \kk$, where $\kk$ is
  the residue field of $Q$. Tensoring with the Koszul complex $K^Q$
  one gets an exact sequence of $Q$-complexes,
  \begin{equation*}
    \tag{$\ast$}
    0 \lra \tp[Q]{\kk}{K^Q} \xra{\a} K^R \xra{\b} K^S \lra 0\:.
  \end{equation*}
  Let $d$ be the dimension of $Q$. From the sequence in homology
  associated to $(\ast)$ one gets the following exact sequence
  \begin{gather*}
    0 \lra \kk \xra{\H[d]{\a}} \H[d]{K^R} \xra{\H[d]{\b}}\H[d]{K^S}\:.
  \end{gather*}
  The rings $R$ and $S$ are artinian, and a rank count yields
  \begin{equation*}
    \type[](R) \deq \rnk[\kk]{(\H[d]{K^R})} \dle \rnk[\kk]{(\H[d]{K^S})} + 1 
    \deq \type[](S) + 1\:.
  \end{equation*}
  It is clear that the residue classes $\cls{g_0}$ and
  $\cls{s_1},\ldots,\cls{s_t}$ in $R$ are non-zero socle
  elements. Moreover, they are $\kk$-linearly independent: Indeed, the
  elements $\cls{s_1},\ldots,\cls{s_t}$ are $\kk$-linearly
  independent, because of the inclusion $\fa \subset \fg$. Further,
  suppose one has $\cls{g_0} = \sum_{i=1}^t \cls{u_i}\cls{s_i}$ where
  the elements $u_i$ are units in $Q$. It follows that
  $g_0-\sum_{i=1}^tu_is_i$ is in $\fa\subseteq \fg$, and as $g_0\in
  \fg$ one gets $\sum_{i=1}^tu_is_i \in \fg$, a contradiction.  Thus,
  there are $t + 1$ $\kk$-linearly independent elements in the
  socle~of~$R$.
\end{prf*}

For the next result, recall from work of J.~Watanabe \cite{JWt73} that
a grade 3 Gorenstein ideal in a regular ring is minimally generated by
an odd number of elements.

\begin{thm}
  \label{thm:2}
  Let $(Q,\n)$ be a regular local ring of dimension $3$ and let
  \mbox{$\fg \subseteq \n^2$} be an $\n$-primary Gorenstein ideal
  minimally generated by elements $g_0,\ldots,g_{2m}$. The ideal $\fa
  = \n g_0 + (g_1,\ldots,g_{2m})$ is $\n$-primary, one has
  $\type[]{(Q/\fa)} =2$ and:
  \begin{prt}
  \item If $m=1$, then $\mu(\fa)=5$ and $Q/\fa$ is of class $\clB$.
  \item If $m=2$, then one of the following holds:
    \begin{itemize}
    \item $\mu(\fa)=4$ and $Q/\fa$ is of class $\clH{3,2}$.
    \item $\mu(\fa)=5$ and $Q/\fa$ is of class $\clB$.
    \item $\mu(\fa) \in \set{6,7}$ and $Q/\fa$ is of class $\clG{r}$
      with $\mu(\fa)-2\ge r \ge \mu(\fa)-3$.
    \end{itemize}
  \item If $m \ge 3$, then $Q/\fa$ is of class $\clG{r}$ with
    $\mu(\fa)-2\ge r \ge \mu(\fa)-3$.
  \end{prt}
\end{thm}

\begin{prf*}
  As $\fg$ defines a Gorenstein ring, one has $\fg:(\fg :\fb) = \fb$
  for every ideal $\fb$ in $Q$ that contains $\fg$.  Let $s\in Q$ be a
  representative of the socle of $Q/\fg$; in $Q$ one has
  \begin{equation*}
    \fg \subseteq (\fa:\n) \subseteq (\fg:\n) = \fg + (s)\:.
  \end{equation*}
  Forming colon ideals one gets $\fg:(\fa:\n) \supseteq \fg:(\fg:\n) =
  \n$ and hence $\fg:(\fa:\n) = \n$.  Forming colon ideals a second
  time now yields $(\fa:\n) = (\fg:\n) = \fg + (s)$; in particular,
  one has $\n s \subseteq \fa$, so it follows from \lemref{type} that
  $\fa$ is $\n$-primary and $R=Q/\fa$ has type $2$; in particular, $R$
  is not Gorenstein.

  Note that one has
  \begin{equation*}
    2m \le \mu(\fa) \le 2m + 3\:.
  \end{equation*}
  Set $S=Q/\fg$; there is an exact sequence of $Q$-modules
  \begin{equation*}
    0 \lra \fg/\fa \lra R \lra S \lra 0
  \end{equation*}
  and an isomorphism $\fg/\fa\is \kk$. Tensor with the Koszul complex
  $K^Q$ to get an exact sequence of $Q$-complexes,
  \begin{equation*}
    \tag{$\ast$}
    0 \lra \tp[Q]{\kk}{K^Q} \xra{\a} K^R \xra{\b} K^S \to 0\;,
  \end{equation*}
  where $\b$ is a morphism of DG $Q$-algebras. Set $A = \H{K^R}$ and
  $G = \H{K^S}$,~one~has
  \begin{gather*}
    \rnk[\kk]{(G_0)} = 1 = \rnk[\kk]{(G_3)} \qqand
    \rnk[\kk]{(G_1)} = 2m+1 = \rnk[\kk]{(G_2)} \\
    \rnk[\kk]{(A_0)} = 1, \ \rnk[\kk]{(A_1)} = \mu(\fa), \
    \rnk[\kk]{(A_2)} = \mu(\fa) +1, \ \: \text{and} \ \:
    \rnk[\kk]{(A_3)} = 2\:,
  \end{gather*}
  and consider the exact sequence in homology associated to $(\ast)$
  \begin{equation*}
    \xymatrix@C=1.8pc@R=2.2pc{
      0 \ar[r] & \kk \ar[rr]^-{\H[3]{\a}}_-{1} && A_3 \ar[rr]_-{1}^-{\H[3]{\b}} &&
      G_3 
      \ar `r[rd]`[l] `[llllld]^-{0} `[dl] [dllll] \\
      & \kk^3 \ar[rr]_-{3}^-{\H[2]{\a}} && A_2 \ar[rr]_-{\mu(\fa)-2}^-{\H[2]{\b}} && G_2 
      \ar `r[rd] `[l] `[llllld]^-{2m+3-\mu(\fa)} `[dl] [dllll] & \\
      & \kk^3 \ar[rr]_-{\mu(\fa)-2m}^-{\H[1]{\a}} && A_1 \ar[rr]_{2m}^-{\H[1]{\b}} && G_1 
      \ar `r[rd] `[l] `[llllld]^-{1} `[dl] [dllll] & \\
      & \kk \ar[rr]_-{0}^-{\H[0]{\a}} && A_0 \ar[rr]_-{1}^-{\H[0]{\b}} && G_0 \ar[r] & 0}
  \end{equation*}
  where the numbers below the arrows indicate the ranks of the
  maps. As $\H{\b}$ is a homomorphism of graded $\kk$-algebras, there
  is a commutative diagram
  \begin{equation*}
    \xymatrix@=3pc{
      G_2 \ar[r]^-{\delta_G}_-{2m+1} & \Hom[\kk]{G_1}{G_3} 
      \ar[d]_-{2m}^-{\Hom[\kk]{\H[1]{\b}}{G_3}} \\
      A_2 \ar[u]_{\mu(\fa)-2}^-{\H[2]{\b}} \ar[r]^-{\varepsilon} 
      & \Hom[\kk]{A_1}{G_3} \\
      A_2 \ar@{=}[u] \ar[r]^-{\delta_A}_-{r} 
      & \ar@{->>}[u]_-{\Hom[\kk]{A_1}{\H[3]{\b}}} \Hom[\kk]{A_1}{A_3} \\
    }
  \end{equation*}
  The rank of $\delta_A$ is at least the rank of $\varepsilon =
  \Hom[\kk]{A_1}{\H[3]{\b}}\circ \delta_A$. The rank of $\varepsilon$
  is at least $(\mu(\fa)-2)-1$ as $\delta_G$ is an isomorphism, see
  \cite[1.4.2]{LLA12}, and the kernel of $\Hom[\kk]{\H[1]{\b}}{G_3}$
  has rank $1$. Thus, one has $r = \rnk[]{(\delta_A)} \ge \mu(\fa)-3$.

  In case $\mu(\fa) \ge 6$, one has $r \ge 3$, and since the type of
  $R$ is $2$, this implies that $R$ is of class $\clG{r}$; see
  \eqref{tab}. This proves part (c) and the last case of part (b).
  For $m=2$ one has $4\le \mu(\fa) \le 7$. If $\mu(\fa)=5$ it follows
  from \prpref{mu5} and \eqref{tab} that $R$ is of class $\clB$. If
  $\mu(\fa) = 4$, then $R$ is of class $\clH{3,2}$ by
  \cite[3.4.2.(a)]{LLA12}. Finally, part (a) is a result of Faucett
  \cite{JAF}.
\end{prf*}

\section{A family of graded local rings of class $\clG{r}$}

\noindent
A grade 3 Gorenstein ideal of a local ring is by a result of Buchsbaum
and Eisenbud \thmcite[2.1]{DABDEs77} minimally generated by
the 
sub-maximal Pfaffians of a $(2m+1) \times (2m+1)$ skew-symmetric
matrix. Thus, skew-symmetric matrices are a source of Gorenstein rings
and, via \thmref{2}, also a source of rings of class $\clG{r}$ that
are not Gorenstein. In this section, we construct an infinite family
of such rings.

\begin{ipg}
  \label{stp:Pfaffian ideal}
  Let $\kk$ be a field and set $Q=\pows[\kk]{x,y,z}$; let $m$ be a
  positive integer.

  Denote by $U_m$ the $m\times m$ matrix over $Q$ whose $i^{\rm th}$
  row has entries
  \begin{equation*}
    u_{i,m-i}=x\,,\quad u_{i,m-i+1}=z\,, \ \text{ and }\ \quad u_{i,m-i+2}=y\,
  \end{equation*}
  and $0$ elsewhere; set
  \begin{equation*}
    d_{-1}=0\,,\quad d_0=1\,,\quad\text{and}\quad d_m = \det (U_m)\:.
  \end{equation*}
  That is,
  \begin{alignat*}{4}
    U_1 &= [z]\;,& \quad 
    U_2 &=
    \begin{bmatrix}
      x & z \\
      z & y
    \end{bmatrix},& \quad U_3 &=
    \begin{bmatrix}
      0 & x & z \\
      x & z & y \\
      z & y & 0
    \end{bmatrix},& \quad U_4 &=
    \begin{bmatrix}
      0 & 0 & x & z \\
      0 & x & z & y\\
      x & z & y & 0\\
      z & y & 0 & 0
    \end{bmatrix}, \quad \ldots\\[2ex]
    d_1 &= z\:,& d_2 &= xy-z^2\:,& d_3 &= 2xyz-z^3\:,& d_4 &= -3xyz^2
    + x^2y^2 + z^4\:, \quad \ldots
  \end{alignat*}
  Notice that for every $i$ in the range $2,\ldots,m$ one has,
  \begin{equation}
    \label{eq:Um}
    U_{m} \deq \left[
      \begin{array}[c]{c|c}
        \Ox & U_{{i-1}}\\[1pt]
        \hline
        \phantom{\raisebox{3pt}{$|$}}U_{m-i+1}\phantom{\raisebox{3pt}{$|$}} & \Oy
      \end{array}
    \right],
  \end{equation}
  where $\Ox$ is the appropriately sized matrix with $x$ in the lower
  right corner and $0$ elsewhere, and $\Oy$ is the matrix with $y$ in
  the top left corner and $0$ elsewhere.

  Let $V_m$ be the $(2m+1) \times (2m+1)$ skew-symmetric matrix given
  by
  \begin{equation}
    \label{eq:Vmd}
    V_m \deq \left[
      \begin{array}[c]{c|c|c}
        O & \Ox & U_m\\[1pt]
        \hline
        -(\Ox)^T & \mspace{9mu}\phantom{\raisebox{3pt}{$|$}}0
        \phantom{\raisebox{3pt}{$|$}}\mspace{9mu} & \Oy \\[1pt]
        \hline
        -U_m & -(\Oy)^T & \mspace{6mu}\phantom{\raisebox{3pt}{$|$}}O
        \phantom{\raisebox{3pt}{$|$}}\mspace{6mu}
      \end{array}
    \right],
  \end{equation}
  where $O$ is the $m\times m$ zero-matrix and, as above, $\Ox$ and
  $\Oy$ are appropriately sized matrices with $0$ everywhere but in
  the lower left and upper right corner, respectively. That is,
  \begin{equation}
    \label{eq:V1}
    V_1 = 
    \begin{bmatrix}
      0 & x & z \\
      -x & 0 & y \\
      -z & -y & 0
    \end{bmatrix}, \quad V_2 =
    \begin{bmatrix}
      0 & 0 & 0 & x & z \\
      0 & 0 & x & z & y \\
      0 & -x & 0 & y & 0\\
      - x & - z & -y & 0 & 0\\
      - z & -y & 0 & 0 & 0
    \end{bmatrix}, \quad \ldots
  \end{equation}

  The sub-maximal Pfaffians of $V_m$ are determined (up to
  a sign) by minors, $\pf{i}{V_m}^2 = \det((V_m)_{ii})$. Consider the
  ideal of $Q$ generated by these Pfaffians,
  \begin{equation}
    \label{eq:pf}
    \fg_m \deq (\pf{1}{V_m},\ldots,\pf{2m+1}{V_m})\:.
  \end{equation}
\end{ipg}

\begin{lem}
  \label{lem:d}
  In the notation from {\rm \pgref{stp:Pfaffian ideal}} the next
  equalities hold for every $m\ge 1$.
  \begin{align*}
    d_m &\deq (-1)^{m-1}z d_{m-1} + xyd_{m-2}\qand\\
    d_m &\deq \sum_{j=0}^{\lfloor\frac{m}{2}\rfloor}
    \binom{m-j}{j}(-1)^{\lfloor\frac{m-2j}{2}\rfloor}x^jy^jz^{m-2j}\:.
  \end{align*}
\end{lem}

\begin{prf*}
  Per \eqref{Um} with $i=2$, expansion of the determinant of $U_m$
  along the first row yields
  \begin{equation*}
    d_m = (-1)^mx\det((U_m)_{1,m-1}) + (-1)^{m+1}z\det (U_{m-1})\:.
  \end{equation*}
  From \eqref{Um} with $i=3$ it follows that expansion along the last
  column yields
  \begin{equation*}
    \det((U_m)_{1,m-1}) = (-1)^my\det (U_{m-2})\:.
  \end{equation*}
  Combining these two expressions, one gets the first equality. The
  second equality now follows by induction.
\end{prf*}

Evidently, the ideal $\fg_m$ from \eqref{pf} is contained in $\n^m$;
in fact, one has $\fg_1 = \n$. One can check that, though the
generating matrices are different, the family of ideals
$\set{\fg_m}_{m\ge 2}$ is the same as that provided by
\prpcite[6.2]{DABDEs77}. To understand what happens when one trims
these ideals, we provide a more detailed description.

\begin{prp}
  \label{prp:dPfaffian ideal}
  Adopt the notation from {\rm \pgref{stp:Pfaffian ideal}} and let
  $\n$ denote the maximal ideal of $Q$. For every $m\ge 2$ the ideal
  $\fg_m \subseteq \n^2$ is an $\n$-primary Gorenstein ideal minimally
  generated by the elements
  \begin{equation*}
    x^{m-i}d_i \ \text{ and } \ y^{m-i}d_i \ \text{ for } \ 
    0 \le i \le m-1\qand d_m\:.
  \end{equation*}
  The ring $Q/\fg_m$ has socle generated by the class of
  $x^{m-1}y^{m-1}$ and Hilbert series
  \begin{equation*}
    \Hilb{Q/\fg_m} = \sum_{i=0}^{m-2} \binom{i+2}{2}\left (t^i+t^{2m-2-i}\right) 
    + \binom {m+1}{2}t^{m-1}\:.
  \end{equation*}
\end{prp}

\begin{prf*}
  Per \eqref{V1} the Pfaffians of $V_1$ are, up to signs,
  \begin{equation*}
    \pf{1}{V_1} = y = yd_0\,, 
    \quad \pf{2}{V_1} = z = d_1\,,\qand\,\pf{3}{V_1} = x =
    xd_0\:.
  \end{equation*}
  For $m \ge 2$ we argue that, up to signs, one has
  \begin{align*}
    \pf{i}{V_m} &= y^{m-i+1}d_{i-1} \ \text{ for } 1\le i \le m\:,\\
    \pf{m+1}{V_m} &= d_m\:,\text{ and}\\
    \pf{2m+2-i}{V_m} &= x^{m-i+1}d_{i-1} \ \text{ for } 1\le i \le m\:.
  \end{align*}
  First notice that the equality $\pf{m+1}{V_m} = d_m$ is immediate
  from \eqref{Vmd}. Further, note that by symmetry in $x$ and $y$ it
  is sufficient to prove that $\pf{i}{V_m} = y^{m-i+1}d_{i-1}$ holds
  for $1\le i \le m$. To compute $\pf{1}{V_m}$ notice that the matrix
  $(V_m)_{11}$ is a $2m\times 2m$-matrix with $\pm y$ on the
  anti-diagonal and zeros below it. Thus, one has $\pf{1}{V_m} = y^m =
  y^md_0$. Now, for $i$ in the range $2,\ldots,m$ consider the matrix
  $(V_m)_{ii}$ as $2\times 2$ block matrix with blocks of size
  $m\times m$,
  \begin{equation*}
    (V_{m})_{ii} \deq \left[
      \begin{array}[c]{c|c}
        X & W_{i}\\[1pt]
        \hline
        \phantom{\raisebox{3pt}{$|$}}-W_i^T\phantom{\raisebox{3pt}{$|$}} & O
      \end{array}
    \right],
  \end{equation*}
  where $O$ is as in \eqref{Vmd}, i.e.\ it is zero. Thus, one has
  \begin{equation*}
    \det((V_{m})_{ii}) \deq 
    \left|
      \begin{array}[c]{c|c}
        X & W_{i}\\[1pt]
        \hline
        \phantom{\raisebox{3pt}{$|$}}-W_i^T\phantom{\raisebox{3pt}{$|$}} & O
      \end{array}
    \right| 
    \deq
    (-1)^m\left|
      \begin{array}[c]{c|c}
        W_i & X\\[1pt]
        \hline
        \phantom{\raisebox{3pt}{$|$}}O\phantom{\raisebox{3pt}{$|$}} & -W_i^T
      \end{array}
    \right| 
    \deq (\det (W_i))^2\:.
  \end{equation*}
  Next, notice that $W_i$ is obtained from $U_m$ by removing row $i$
  and adding a row $\Oy$ at the bottom. Thus, per \eqref{Um} it has
  the form
  \begin{equation*}
    W_i \deq \left[
      \begin{array}[c]{c|c}
        \Ox & U_{i-1}\\[1pt]
        \hline
        \phantom{\raisebox{3pt}{$|$}}Y\phantom{\raisebox{3pt}{$|$}} & O
      \end{array}
    \right],
  \end{equation*}
  where $Y$ is the matrix obtained from $U_{m-i+1}$ by removing the
  first row and adding a row $\Oy$ at the bottom. In particular, it is
  a $(m-i+1)\times (m-i+1)$-matrix with $\pm y$ on the anti-diagonal
  and zeros below it.  Thus, computing the determinant of $W_i$ by
  successive expansion on the last $m-i+1$ rows one gets, up to a
  sign, $\pf{i}{V_m} = y^{m-i+1}d_{i-1}$. It follows that $\fg_m$ is
  generated by the listed elements.

  The elements $x^m$, $y^m$, $d_m$ form a $Q$-regular sequence in
  $\fg_m$, so it follows from \thmcite[2.1]{DABDEs77} that $\fg_m$ is
  a Gorenstein ideal minimally generated by the listed elements. In
  particular, $\fg_m$ is $\n$-primary. In fact, in this case it is
  elementary to see that the generating set is minimal: Notice from
  \lemref{d} that $d_i$ is a linear combination of monomials of the
  form $x^jy^jz^{i-2j}$. Hence, each generator $x^{m-i}d_i$ is a
  linear combination of monomials of the form $x^{m-i+j}y^jz^{i-2j}$
  while the generators $y^{m-i}d_i$ are linear combinations of
  monomials $x^jy^{m-i+j}z^{i-2j}$. Thus the generators are linear
  combinations of disjoint sets of degree $m$ monomials and hence
  linearly independent.

  The Hilbert series of the power series ring $Q$ is
  $\Hilb{Q}=\sum_{i=0}^{\infty}\binom{i+2}{2}t^i$.  Since $\fg_m$ is
  Gorenstein and minimally generated by $2m+1$ elements of degree $m$,
  the Hilbert series of the ring $S_m=Q/\fg_m$ is symmetric and given
  by
  \begin{equation*}
    \Hilb{S_m} = \sum_{i=0}^{m-2} \binom{i+2}{2}\left (t^i+t^{2m-2-i}\right) 
    + \binom {m+1}{2}t^{m-1}\:.
  \end{equation*}
  In particular, the socle degree of $S_m$ is $2m-2$. Evidently, one
  has $(x^{m-1}y^{m-1})\mathfrak{n} \subseteq \fg_m$, so it is
  sufficient to show that the element $x^{m-1}y^{m-1}$ is not in
  $\fg_m$, i.e.\ that it yields a non-zero socle element in $S_m$. If
  it were in $\fg_m$, then one would have $x(x^{m-2}y^{m-1})$ in
  $\fg_m$ along with $x^{m-2}(y^md_0) = y(x^{m-2}y^{m-1})$ and
  $x^{m-2}(y^{m-1}d_1) = z(x^{m-2}y^{m-1})$. Thus, $x^{m-2}y^{m-1}$
  would yield a socle element in $S_m$ of degree $2m-3$, whence it
  must be $0$; i.e.\ one would have $x^{m-2}y^{m-1} \in
  \fg_m$. Reiterating this argument, one arrives at the conclusion
  that $y^{m-1}$ is in $\fg_m$, which is absurd as the generators of
  $\fg_m$ have degree $m$.
\end{prf*}

Finally, we apply the trimming procedure from \thmref{2} to the ideals
$\fg_m$.  \enlargethispage*{\baselineskip}
\begin{ipg}
  \label{rmk:a}
  Adopt the notation from {\rm \pgref{stp:Pfaffian ideal}}. By
  \prpref{dPfaffian ideal} one has
  \begin{equation*}
    \fg_2 \deq (x^2, xz, xy-z^2, yz, y^2)\;.
  \end{equation*}
  Trimming the generators $xz$ and $yz$ one gets the following ideals
  of $Q$,
  \begin{align*}
    (x,y,z)xz + (x^2, xy-z^2, yz, y^2) &\deq (x^2, xy-z^2, yz, y^2)\quad\text{and}\\
    (x,y,z)yz + (x^2, xz, xy-z^2, y^2) &\deq (x^2, xz, xy-z^2, y^2)\:.
  \end{align*}
  They are both minimally generated by 4 elements, so they define
  quotient rings of class $\clH{3,2}$; see \thmref{2}(b).  Moreover,
  one has
  \begin{align*}
    (x,y,z)x^2 + (xz, xy-z^2, yz, y^2) &\deq (x^3, xz, xy-z^2, yz, y^2)\:,\\
    (x,y,z)y^2 + (x^2,xz, xy-z^2, yz) &\deq (x^2, xz, xy-z^2, yz,
    y^3)\:,
    \quad\text{and}\\
    (x,y,z)(xy-z^2) + (x^2, xz, yz, y^2) &\deq (x^2, xz, z^3, yz,
    y^2)\:,
  \end{align*}
  so by \thmref{2}(b) these ideals define rings of class $\clB$.
\end{ipg}

From the next result one immediately gets the statement of \thmref{a}
about existence of infinite families of rings of class $\clG{r}$ that
are not Gorenstein.

\begin{prp}
  \label{prp:a}
  Adopt the notation from {\rm \pgref{stp:Pfaffian ideal}} and let
  $\n$ denote the maximal ideal of $Q$.  Let $g$ be one of the
  generators of $\fg_m$ listed in {\rm \pgref{prp:dPfaffian ideal}},
  let $\fb$ be the ideal generated by the remaining $2m$ generators of
  $\fg_m$, and set $\fa = \n g + \fb$.  For $m \ge 3$ the ring
  $R=Q/\fa$ has the following properties.
  \begin{prt}
  \item $R$ is an artinian local ring of type 2 with socle generated
    by the classes of the elements $g$ and $x^{m-1}y^{m-1}$.
  \item If $g$ is $x^{m-i}d_i$ or $y^{m-i}d_i$ for some
    $i\in\set{1,\ldots,m-1}$, then $\fa$ is minimally generated by
    $2m$ elements and $R$ is of class $\clG{2m-3}$.
  \item If $g$ is $x^m$, $y^m$, or $d_m$, then $\fa$ is minimally
    generated by $2m+1$ elements and $R$ is of class $\clG{2m-2}$.
  \end{prt}
\end{prp}

\begin{prf*}
  Fix $m\ge 3$; for brevity the class in $R$ or $S=Q/\fg_{m}$ of an
  element $u$ in $Q$ is also written $u$.

  Part (a) is immediate from \lemref{type}.  We prove parts (b) and
  (c) together. First we describe the generators of $\fa$ using the
  recurrence formula from \lemref{d}. For $1 \le i \le m$ one has
  \begin{equation}
    \label{eq:1}
    \begin{aligned}
      x(x^{m-i}d_i) & = x^{m-(i-1)}((-1)^{i-1}zd_{i-1} + xyd_{i-2}) \\
      & = (-1)^{i-1}z(x^{m-(i-1)}d_{i-1}) + y(x^{m-(i-2)}d_{i-2})\:.
    \end{aligned}
  \end{equation}
  For $0 \le i \le m-2$ one has
  \begin{equation}
    \label{eq:2}
    \begin{aligned}
      y(x^{m-i}d_i) & = x^{m-(i+1)}(xyd_i) \\
      &= x^{m-(i+1)}(d_{i+2} - (-1)^{i+1}zd_{i+1}) \\
      &= x(x^{m-(i+2)}d_{i+2}) + (-1)^{i}z(x^{m-(i+1)}d_{i+1}) \qtext{and moreover}\\
      y(xd_{m-1}) & = x(yd_{m-1})\:.
    \end{aligned}
  \end{equation}
  For $0 \le i \le m-1$ one has
  \begin{equation}
    \label{eq:3}
    \begin{aligned}
      z(x^{m-i}d_i) & = x^{m-i}(-1)^i(d_{i+1}-xyd_{i-1}) \\
      &= (-1)^ix(x^{m-(i+1)}d_{i+1}) - (-1)^iy(x^{m-(i-1)}d_{i-1})\:.
    \end{aligned}
  \end{equation}

  For $g=x^{m-i}d_i$ with $1\le i \le m-1$ it follows immediately from
  \eqref{1}--\eqref{3} that $\n g$ is contained in $\fb$ , so $\fa =
  \fb$ is minimally generated by $2m$ elements. By symmetry the same
  is true for $g=y^{m-i}d_i$ with $1\le i \le m-1$.

  For $g=x^m$ one has $yg \in \fb$ and $zg \in \fb$ by \eqref{2} and
  \eqref{3}, so $\fa$ is generated by the $2m$ generators of $\fb$ and
  $x^{m+1}$. To see that this is a minimal set of generators, note
  that the generators of $\fb$ have degree $m$ and none of them
  includes the term $x^m$. The statement for $g=y^m$ follows by
  symmetry.

  For $g=d_m$ one has $xg \in \fb$ by \eqref{1} and $yg \in \fb$ by
  symmetry.  Thus $\fa$ is generated by the $2m$ generators of $\fb$
  and $zd_m$. To see that this is a minimal set of generators, note
  from \lemref{d} that $zd_m$ has a $z^{m+1}$ term, while the
  generators of $\fb$ have degree $m$ and none of them has a $z^m$
  term.

  To determine the multiplicative structure on $A = \H[]{K^R}$ we
  first describe a basis for $A_1$. The Koszul complex $K^R$ is the
  exterior algebra of the free $R$-module with basis $\set{\e{x},
    \e{y}, \e{z}}$ endowed with the differential given by
  $\partial(\e{x}) = x$, $\partial(\e{y}) = y$, and $\partial(\e{z}) =
  z$. We suppress the wedge in products on $K^R$ and adopt the
  following shorthands
  \begin{equation*}
    \e{xy} = \e{x}\e{y}\,,\quad
    \e{xz} = \e{x}\e{z}\,, \quad
    \e{yz} = \e{y}\e{z}\,,\qand\, \e{xyz} = \e{x}\e{y}\e{z}\,.
  \end{equation*}
  Because of the symmetry in $x$ and $y$ we only consider
  $g=x^{m-i}d_i$. Given the minimal generating set of $\fa$ described
  above, one gets:

  If $g=x^m$ then the following cycles in $K^R_1$ yield a basis for
  $A_1$
  \begin{align*}
    x^m\e{x} \qand x^{m-j-1}d_j\e{x} & \text{ for } 1 \le j \le m-1, \\
    y^{m-j-1}d_j\e{y} & \text{ for } 0 \le j \le m-1,\qand\\
    (-1)^{m-1}z^{m-1}d_{m-1}\e{z} + xd_{m-2}\e{y}\:.
  \end{align*}
  If $g=x^{m-i}d_i$ for some $i$ in the range $1,\ldots,m-1$, then the
  following cycles in $K^R_1$ yield a basis for $A_1$
  \begin{align*}
    x^{m-j-1}d_j\e{x} & \text{ for } 0 \le j \le m-1,\, j\ne i\\
    y^{m-j-1}d_j\e{y} & \text{ for } 0 \le j \le m-1, \qand\\
    (-1)^{m-1}z^{m-1}d_{m-1}\e{z} + xd_{m-2}\e{y}\:.
  \end{align*}
  If $g=d_m$ then the following cycles in $K^R_1$ yield a basis for
  $A_1$
  \begin{align*}
    x^{m-j-1}d_j\e{x} & \text{ for } 0 \le j \le m-1,\\
    y^{m-j-1}d_j\e{y} & \text{ for } 0 \le j \le m-1,\qand\\
    d_m\e{z}\:.
  \end{align*}
  From \thmref{2} it is known that $R$ is of class $\clG{r}$ with
  $\mu(\fa)-3 \le r$. To prove that equality holds, which is the claim
  in (b) and (c), it suffices to show that the kernel of $\delta$ has
  rank at least $(\mu(\fa) + 1) - (\mu(\fa) - 3) = 4$; see
  \pgref{class}. To this end we first notice that the cycles
  $g\e{xy}$, $g\e{xz}$, and $g\e{yz}$ yield linearly independent
  elements of $A_2$. Assume towards a contradiction that they are not,
  then there exists an element $h\e{xyz}$ in $K_3^Q$ and elements
  $q_1$, $q_2$, and $q_3$ in $Q$ and not all in $\n$ with
  \begin{equation*}
    \partial(h\e{xyz}) - (q_1g\e{xy} + q_2g\e{xz} + q_3g\e{yz}) \in \fa K_2^Q.
  \end{equation*}
  That is, one has $zh -q_1g\in \fa$, $yh + q_2g\in \fa$, and $xh
  -q_3g\in \fa$, and hence $h\not\in\n^m$ as $g \not\in \fa +
  \n^{m+1}$. Furthermore, the class of $h$ is a socle element in $S$
  as one has $\n h \subseteq \fa + Qg = \fg_m$. Thus, $h\in \fg_m$ or
  $h = qx^{m-1}y^{m-1}$ for some $q\in Q\setminus\n$. In either case
  one has $h\in \n^m$, which is a contradiction. Thus $g\e{xy}$,
  $g\e{xz}$, and $g\e{yz}$ yield linearly independent elements in
  $A_2$ that clearly belong to the kernel of $\delta$.

  Finally we produce a fourth element in the kernel.  For $g=x^n$ the
  element
  \begin{equation*}
    f=y^{m-1}\e{yz}
  \end{equation*}
  is clearly a cycle in $K^R_2$, and it is not a boundary. Indeed, if
  one had $f = \partial(h\e{xyz}) = hx\e{yz} -hy\e{xz} + hz\e{xy}$ for
  some homogeneous element $h\in R$, then it would have degree $m-2$
  and one would have $hy=0=hz$ in $R$, which is impossible as $\fa$
  has generators of degree at least $m$. The products
  $(y^{m-j-1}d_j\e{y})\cdot f$ and $((-1)^{m-1}z^{m-1}\e{z} +
  xd_{m-2}\e{y})\cdot f$ in $K^R$ vanish by graded
  commutativity. Moreover, one has
  \begin{align*}
    (x^m\e{x})\cdot f &= x(x^{m-1}y^{m-1})\e{xyz} = 0\qand\\
    (x^{m-j-1}d_j\e{x})\cdot f &= x^{m-j-1}y^{j-1}(y^{m-j}d_j)\e{xyz}
    = 0\:.
  \end{align*}
  Thus the homology class of $f$ annihilates $A_1$.

  For $g=x^{m-i}d_i$ and $1 \le i \le m-1$ the element
  \begin{equation*}
    f = y^{m-i}d_{i-1}\e{xy} + (-1)^{i-1}y^{m-i-1}d_i\e{yz}
  \end{equation*}
  is a cycle in $K^R_2$; indeed one has
  \begin{align*}
    \partial(f) & = xy^{m-i}d_{i-1}\e{y} - y^{m-(i-1)}d_{i-1}\e{x} +
    (-1)^{i-1}y^{m-i}d_i\e{z} + (-1)^{i}y^{m-i-1}zd_i\e{y}\\
    &= y^{m-i-1}((-1)^{i}zd_i + xyd_{i-1})\e{y}\\
    &= y^{m-(i+1)}d_{i+1}\\
    &=0\;,
  \end{align*}
  where the third equality follows from \lemref{d}.  An argument
  similar to the one above shows that $f$ is not a boundary.  The
  products $(y^{m-j-1}d_j\e{y})\cdot f$ in $K^R$ vanish by graded
  commutativity. Moreover, one has
  \begin{align*}
    (x^{m-j-1}d_j\e{x})\cdot f &=
    (-1)^{i-1}x^{m-j-1}d_jy^{m-i-1}d_i\e{xyz}\:.
  \end{align*}
  If $i > j$ holds, then the element $x^{m-j-1}d_jy^{m-i-1}d_i$ is $0$
  in $R$ because it is divisible by $g$, which is a socle element in
  $R$. If one has $i < j$, then the element $x^{m-j-1}d_jy^{m-i-1}d_i$
  is zero in $R$ because it is divisible in $Q$ by the generator
  $y^{m-j}d_j$ of $\fa$. Finally, one has
  \begin{align*}
    ((-1)^{m-1}z^{m-1}d_{m-1}\e{z} + xd_{m-2}\e{y})\cdot f
    &= (-1)^{m-1}y^{m-i}d_{i-1}z^{m-1}d_{m-1}\e{xyz}\\
    &= (-1)^{m-1}y^{m-i-1}d_{i-1}z^{m-1}(yd_{m-1})\e{xyz}\\
    &= 0
  \end{align*}
  in $K^R$, so the homology class of $f$ annihilates $A_1$.

  For $g=d_m$ the element
  \begin{equation*}
    f = d_{m-1}\e{xy}
  \end{equation*}
  is evidently a cycle in $K^R_2$, and as above it is not a boundary.
  The products $(x^{m-j-1}d_j\e{x})\cdot f$ and
  $(y^{m-j-1}d_j\e{y})\cdot f$ in $K^R$ vanish by graded
  commutativity. Finally one has,
  \begin{equation*}
    (d_m\e{z})\cdot f = d_{m-1}d_m\e{xyz}=0\:,
  \end{equation*}
  as $g=d_m$ is a socle element of $R$.
\end{prf*}

\section*{Acknowledgments}

We thank Parangama Sarkar for alerting us to a flaw in a previous version of \lemref{type}.

\bibliographystyle{amsplain}

\providecommand{\MR}[1]{\mbox{\href{http://www.ams.org/mathscinet-getitem?mr=#1}{#1}}}
   \renewcommand{\MR}[1]{\mbox{\href{http://www.ams.org/mathscinet-getitem?mr=#1}{#1}}}
\providecommand{\href}[2]{#2}

\end{document}